\documentclass[12pt,a4paper,reqno]{amsart}
\usepackage{amsfonts}
\usepackage{amssymb}
\usepackage{amsthm}
\usepackage{amsmath}
\usepackage{newlfont}
\usepackage{graphicx}
\usepackage{color}

\def\r{\mathbb R}
 
\def\s{\mathbb S}

\newtheorem{theorem}{Theorem}[section]
\newtheorem{definition}[theorem]{Definition}
\newtheorem{proposition}[theorem]{Proposition}
\newtheorem{remark}[theorem]{Remark}
\newtheorem{corollary}[theorem]{Corollary}
 
\newtheorem{lemma}[theorem]{Lemma} 

\title{Separable surfaces that are critical points of the Dirichlet energy}

\author{Rafael L\'opez}
\address{Departamento de Geometr\'{\i}a y Topolog\'{\i}a\\
 Universidad de Granada\\
 18071 Granada, Spain }
\email{rcamino@ugr.es}
 

\date{}

\begin{document}
\begin{abstract}
In this paper, we study surfaces $z=\varphi(x,y)$ in Euclidean space  that satisfy the equation $\varphi_{xx}+\varphi_{yy}=\frac{\Lambda}{2}$ where $\Lambda\in\r$ is a real constant. We classify these  surfaces when they are the zero level sets  of an implicit equation of the type $f(x)+g(y)+h(z)=0$, where $f$, $g$ and $h$ are smooth functions of one variable.  If $\Lambda=0$, we find a large family of surfaces with interesting  symmetry properties. However, if $\Lambda\not=0$, we show that the surfaces must be either surfaces of revolution or   of the type $z=f(x)+g(y)$; furthermore, explicit parametrizations of these surfaces are obtained.
\end{abstract}

\maketitle

 %%%%%%%%%%%%%%%%%%%%%%%%%%%%%%%%
\section{Introduction and  classification} \label{intro}
%%%%%%%%%%%%%%%%%%%%%%%%%%%%%%%%%%%%%%%%%%%%%%%%%%%%

In this paper, we study solutions of the equation 
\begin{equation}\label{eq1}
\varphi_{xx}+\varphi_{yy}=\frac{\Lambda}{2},
\end{equation}
 where $\varphi=\varphi(x,y)$ is a smooth function defined in a domain  $\Omega\subset\r^2$ and $\Lambda\in\r$ is a real constant. Solutions $\varphi$ of Eq. \eqref{eq1} are critical points of the functional 
$$E[\varphi]=\int_\Omega|D\varphi|^2+\Lambda\int_\Omega \varphi.$$
The first term is the Dirichlet energy of $\varphi$ and the second integral is the volume of the graph   $z=\varphi(x,y)$. Let $\Sigma$ be the  graph $z=\varphi(x,y)$ viewed as a surface in Euclidean space $\r^3$. If  $\nu =(\nu_1,\nu_2,\nu_3)$  is the unit normal vector   to $\Sigma$, then   $\nu=(-D\varphi,1)/\sqrt{1+|D\varphi|^2}$ and thus $\nu_3=1/\sqrt{1+|D\varphi|^2}$ This allows us to express the Dirichlet energy as a functional $\mathcal{F}$ defined by 
\begin{equation}\label{fu}
\mathcal{F}(\Sigma):=\int_\Sigma \left(\dfrac{1}{\nu_3}-\nu_3\right)\, d\Sigma,
\end{equation}
 where $d\Sigma$ is the area element of $\Sigma$. Consequently,  $\Sigma$ is a critical point of  $E$  if and only if $\Sigma$ is a critical point of the energy \eqref{fu} for all   compactly supported, volume-preserving variations of $\Sigma$, where $\Lambda$ is a Lagrange multiplier.  The  functional $\mathcal{F}$ defined in \eqref{fu} is of anisotropic type because the   integral  depends on the unit normal vector $\nu$ of $\Sigma$. In general, anisotropic  functionals  are energies of the type
\begin{equation}\label{fu2}
\mathcal{F}(\Sigma)=\int_\Sigma F(\nu)\, d\Sigma,
\end{equation}
 where $F\colon U\subset\s^2\to\r^+$ is a smooth function defined on an open set $U$ of the sphere $\s^2$. These energies  appear in fluid phenomena when the surface tensions of interfaces  depend  on the orientation \cite{ta}. From a differential-geometric viewpoint, anisotropic  functionals \eqref{fu2} have received significant interest, especially regarding problems related to the stability of the surface. We refer the reader to the works of Koiso and Palmer  \cite{kp2,kp3,kp5} and the references therein.   
 On the other hand, very recently the author has obtained a full classification of the surfaces that are solutions of \eqref{eq1} of ruled type \cite{lo1} or of rotational and cyclic type \cite{lo2}.
 
 For functionals of type \eqref{fu2}, a surface $\Sigma$    is a critical point  for all volume-preserving variations of $\Sigma$ if and only if the quantity $\Lambda$, defined by
$$\Lambda:= 2HF-\mbox{div}_\Sigma DF,$$
 is constant,   where $DF$ is the gradient on $\s^2$ and $H$ is the mean curvature of $\Sigma$.    If the functional \eqref{fu2} is elliptic, we associate it  with the {\it Wulff shape}   $\xi(\nu)=DF(\nu)+F(\nu)\nu$, which is a convex surface. The Wulff shape is the unique global minimizer of $\mathcal{F}$ when a volume constraint is fixed \cite{br,ta}. For the Dirichlet energy   \eqref{fu}   the anisotropic mean curvature coincides with the value $\Lambda$ in Eq. \eqref{eq1} and   the Wulff shape is  the paraboloid    $z=x^2+y^2$.  On the other hand, the isotropic case corresponds to $F\equiv 1$ in \eqref{fu2} where $\Lambda$  is now the mean curvature $H$ of the surface and the Wulff shape is the round sphere $\s^2\subset\r^3$.
 
An approach to finding solutions of Eq.  \eqref{eq1} is to consider a separation of variables of the type   $\varphi(x,y)= f(x)+g(y)$ or $\varphi(x,y)=f(x)g(y)$, where $f$ and $g$ are smooth functions of one variable. These solutions are easily obtained from \eqref{eq1}: see also Prop. \ref{pr1}. 

In this paper, we   address the problem of finding   the solutions of \eqref{eq1} by separation of variables in full generality.  For this purpose, we introduce the following definition.

\begin{definition} A  (regular) surface $\Sigma$ in $\r^3$ is said to be {\it separable} if $\Sigma$  can be expressed as  
\begin{equation}\label{s1}
\Sigma=\{(x,y,z)\in\r^3:f(x)+g(y)+h(z)=0\},
\end{equation}
 where $f$, $g$ and $h$ are   smooth functions defined in   intervals   of $\r$, respectively. 
  \end{definition}
  
 By the regularity of $\Sigma$,    $f'(x)^2+g'(y)^2+h'(z)^2\not=0$ for every $(x,y,z)\in\Sigma$. The symbol $(')$ denotes the derivative with respect to the corresponding variable of the function. In such a case, the unit normal vector is $\nu=(f',g',z')/\sqrt{f'^2+g'^2+h'^2}$. 

For particular choices of functions in \eqref{s1}, we obtain special types of separable surfaces. 
 \begin{enumerate}

 \item  Translation surfaces. A translation surface is a surface given by $z= f(x)+g(y)$. These surfaces correspond to \eqref{s1} when $h(z)=-z$. A translation surface $\Sigma$ is obtained by the sum of the planar curves $x\mapsto (x,0,f(x))$ and $y\mapsto (0,y,g(y))$. A particular case occurs when $f$ or $g$ is a constant function, in which case, $\Sigma$ is a cylindrical  surface; that is, a ruled surface whose rulings are parallel straight lines. For example, if $g(y)=y_0$, then $\Sigma$ is   generated by the planar curve $z=f(x)$ contained in the $xz$-coordinate plane and the rulings are  parallel to the $y$-axis.
 
 \item Homothetical surfaces. These surfaces are given by   $z=p(x)q(y)$. Here, the choice in \eqref{s1} is  $h(z)=-\log z$.
 \item   Rotational surfaces. If $\Sigma$ is a rotational surface    with  respect to   the $z$-axis, then $\Sigma$ can be expressed as  $z=p(x^2+y^2)$. Here $f(x)=x^2$ and $g(y)=y^2$.  
\end{enumerate}

Notation convention:  To increase readability of the solutions obtained in this paper, we use $a_i,b_i,\ldots$ to denote  arbitrary constants of integration arising in the expressions of solutions of ODEs. We only discuss explicitly the cases in which some of these constants are non-zero.

Surfaces satisfying Eq. \eqref{eq1} and being of one of the above three types are obtained by a direct integration. The solutions are described in the following result.

\begin{proposition}\label{pr1}
Let $\Sigma$ be a   surface   satisfying \eqref{eq1}. In the following statements, the roles of the variables $x$ and $y$ can be interchanged.
\begin{enumerate}

\item (Translation surface). If $z=f(x)+g(y)$, then  
\begin{equation}\label{sol2}
\begin{split}f(x)&=\frac{\Lambda+2m}{4}x^2+a_1x+a_2,\\
g(y)&=-\frac{m}{2}y^2+b_1y+b_2,
\end{split}
\end{equation}
with $m\in\r$.
 \item (Homothetical surface). Suppose $z(x,y)=p(x)q(y)$.
 \begin{enumerate}
 \item Case $\Lambda=0$. Then   $p$ and $q$ are linear functions or  
  \begin{equation}\label{sol31}
 \begin{split}p(x)&=a_1\cos(kx)+a_2\sin(kx),\\
q(y)&=b_1\cosh(ky)+b_2\sinh(ky),\\
\end{split}
\end{equation}
with $k\in\r$,  $a_1^2+a_2^2\not=0$ and $b_1^2+b_2^2\not=0$.
 \item Case $\Lambda\not=0$.  Then $q(x)=q_0\not=0 $ is constant and 
 \begin{equation}\label{sol32}
p(x)=\frac{\Lambda}{4q_0}x^2+a_1x+a_2.
\end{equation}

 \end{enumerate}
 
 \item (Rotational surface). If $z=p(x^2+y^2)$, then   
 \begin{equation}\label{sol4}
 z(x,y)= c_1\log( x^2+y^2)+\frac{\Lambda}{8}(x^2+y^2)+c_2.
 \end{equation}
 \end{enumerate}
\end{proposition}

\begin{remark}\label{rem1}
\begin{enumerate}
\item Cylindrical surfaces appear in \eqref{sol2} when $m=b_1=0$, which coincide with the homothetical surfaces given in  \eqref{sol32}.
\item If we take $m=-\Lambda/4$ in the translation surface \eqref{sol2}, then we obtain the rotational surface  \eqref{sol4} for $c_1=0$.
 \item The Wulff shape corresponds to $c_1=0$ in \eqref{sol4}.  
 \item All translation surfaces \eqref{sol2}, and homothetical surfaces \eqref{sol31}, \eqref{sol32} are entire graphs. Rotational surfaces \eqref{sol4} are entire graphs if $c_1=0$ and graphs on $\r^2-\{(0,0)\}$ if $c_1\not=0$. 
 \end{enumerate}
\end{remark}

This work is also motivated by the isotropic setting, corresponding to $F\equiv 1$ in \eqref{fu2}.  The classification of the separable surfaces with constant mean curvature $H$ is well known. When $H=0$, there exists a large variety  of examples due to Schwarz, Fr\'echet, Scherk, and Weingarten, among others: see \cite[pp. 71-80]{ni} for a description of them. When $H\not=0$, the classification was obtained more recently by T. Hasanis and the author, who proved that the only separable surfaces with non-zero constant mean curvature are surfaces of revolution (Delaunay surfaces)  \cite{hl}.

In this paper, we classify all  separable surfaces that are critical points of the Dirichlet energy.

 \begin{theorem}\label{t1}
 Let $\Sigma$ be a separable surface given by the implicit equation \eqref{s1}. If $\Sigma$ satisfies \eqref{eq1}, then $\Sigma$ is a surface described in Prop. \ref{pr1} or $\Lambda=0$ and the  first integrals of the functions $f$, $g$ and $h$ in \eqref{s1}   are given by   
   \begin{equation}\label{ss1}
   \left\{
  \begin{split}
 f'(x)^2&=r\cos (k f(x)-d_1) +a,\\
g'(y)^2&=r\cos (k g(y)-d_2)-a,\\
h'(z)^2&=c (1+\cos( k h(z)+d_1+d_2),
\end{split}\right.
\end{equation} 
or 
   \begin{equation}\label{ss2}
   \left\{
  \begin{split}
 f'(x)^2&=r\cosh (k f(x)-d_1) +a,\\
g'(y)^2&=r\cosh (k g(y)-d_2)-a,\\
h'(z)^2&=c (1+\cosh( k h(z)+d_1+d_2),
\end{split}
\right.
\end{equation} 
with  $r,k,c\not=0$.
 \end{theorem}

This classification exhibits similarities with the isotropic case. If we distinguish whether $\Lambda $ is $0$ or not, we see that if $\Lambda=0$,  the family of separable surfaces is large, with many examples obtained after fixing the constants of integration in  \eqref{ss1} and \eqref{ss2} (see   Sect. \ref{s5}). In contrast, if $\Lambda\not=0$, then Thm. \ref{t1} shows that the only possibilities are those already described in   Prop. \ref{pr1}. Taking into account item (1) of Rem. \ref{rem1}, we emphasize this case in the following corollary.

\begin{corollary} 
If $\Lambda\not=0$, then the surfaces defined in \eqref{sol2} and \eqref{sol4} are the only separable surfaces    satisfying Eq. \eqref{eq1}. 
\end{corollary}

The organization of the paper is as follows. In Sect. \ref{s2} we prove Prop. \ref{pr1}. We prove Thm. \ref{t1} in two steps, which will be carried out in  Sects. \ref{s3} and \ref{s4}. In the last section (Sect. \ref{s5}), we show explicit examples of surfaces obtained via Thm. \ref{t1} by solving  \eqref{ss1} and \eqref{ss2}.    These surfaces have interesting  symmetry properties, proving in particular that they are  doubly periodic surfaces.

%%%%%%%%%%%%
\section{Proof of Proposition \ref{pr1}}\label{s2}
%%%%%%%%%%%%%%%%%%

\begin{figure}[hbtp]
\begin{center}\includegraphics[width=.3\textwidth]{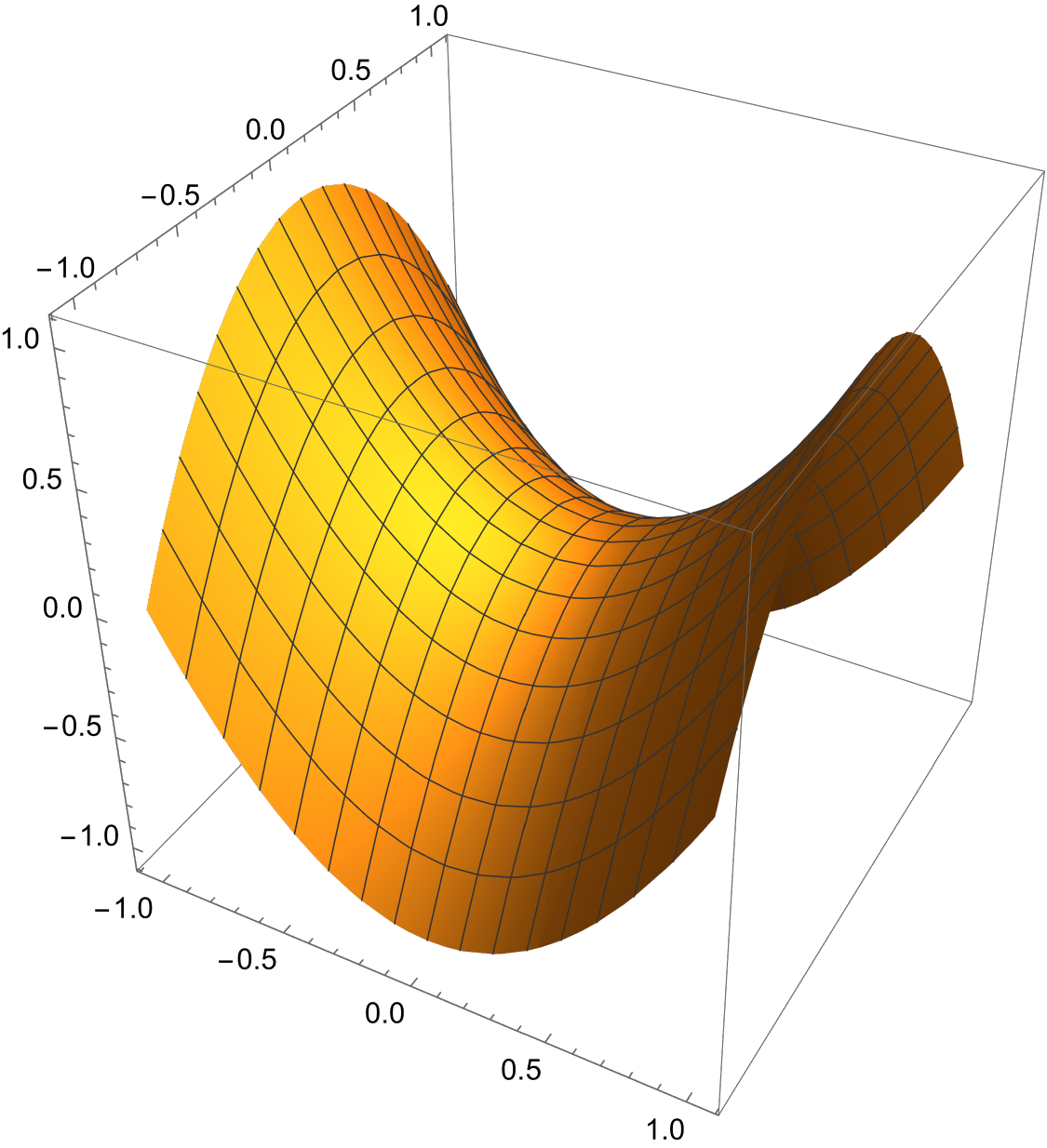}\quad \includegraphics[width=.31\textwidth]{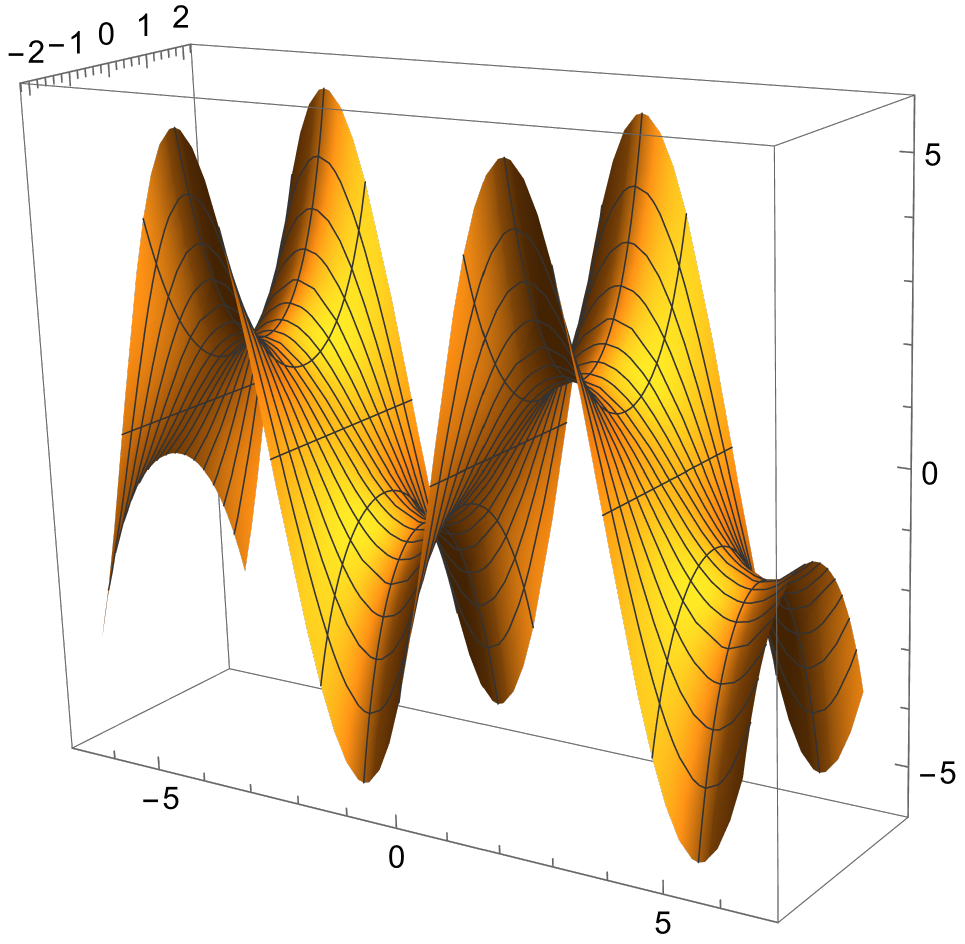}
\quad \includegraphics[width=.31\textwidth]{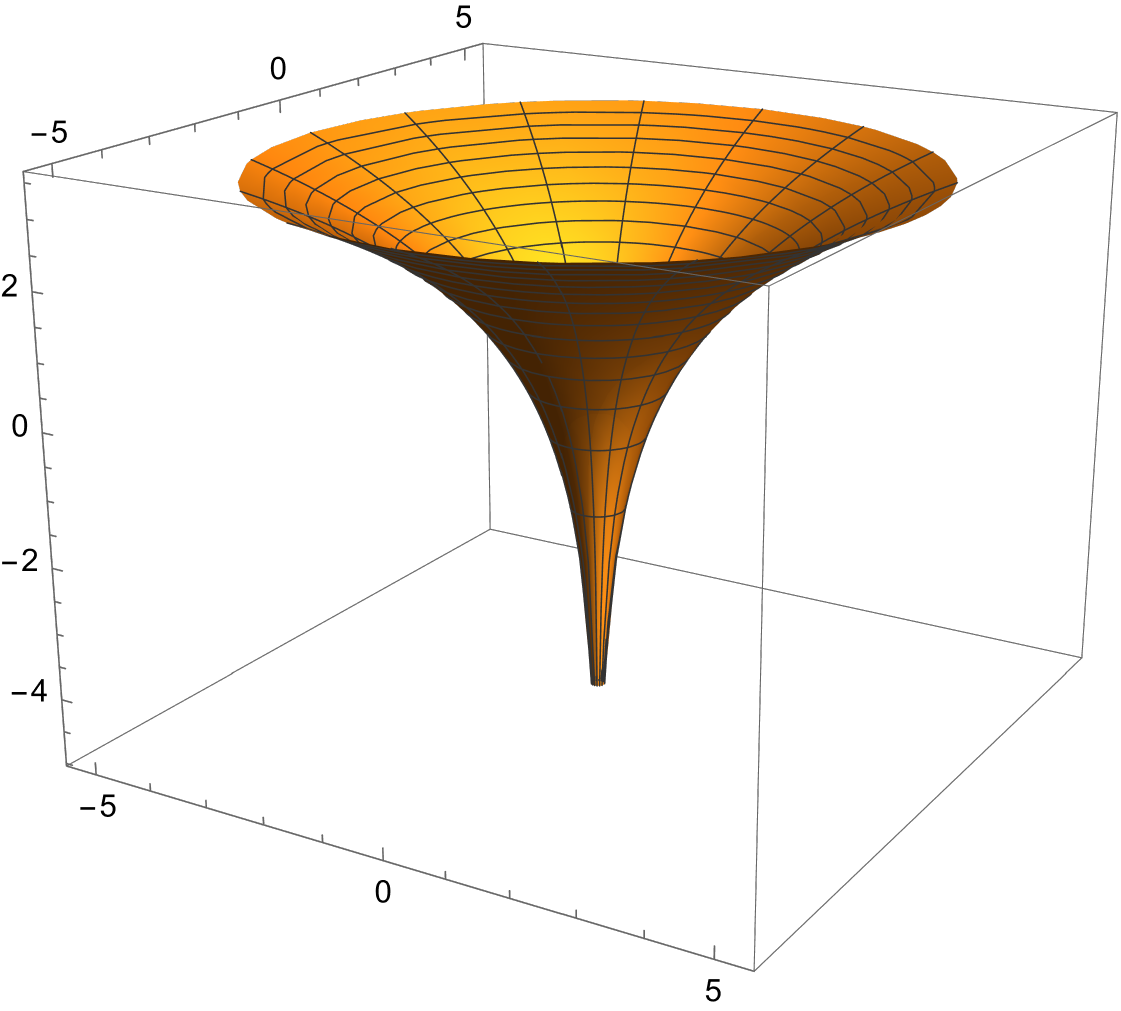} \end{center}
\caption{Separable surfaces of Prop. \ref{pr1} in the case  $\Lambda=0$: translation surface  (left), homothetical surface $z=(\sin x-\cos x)\cosh y$ (middle); and rotational surface $z=\log(x^2+y^2)$ (right). }\label{fig1}
\end{figure}

\begin{figure}[hbtp]
\begin{center}\includegraphics[width=.26\textwidth]{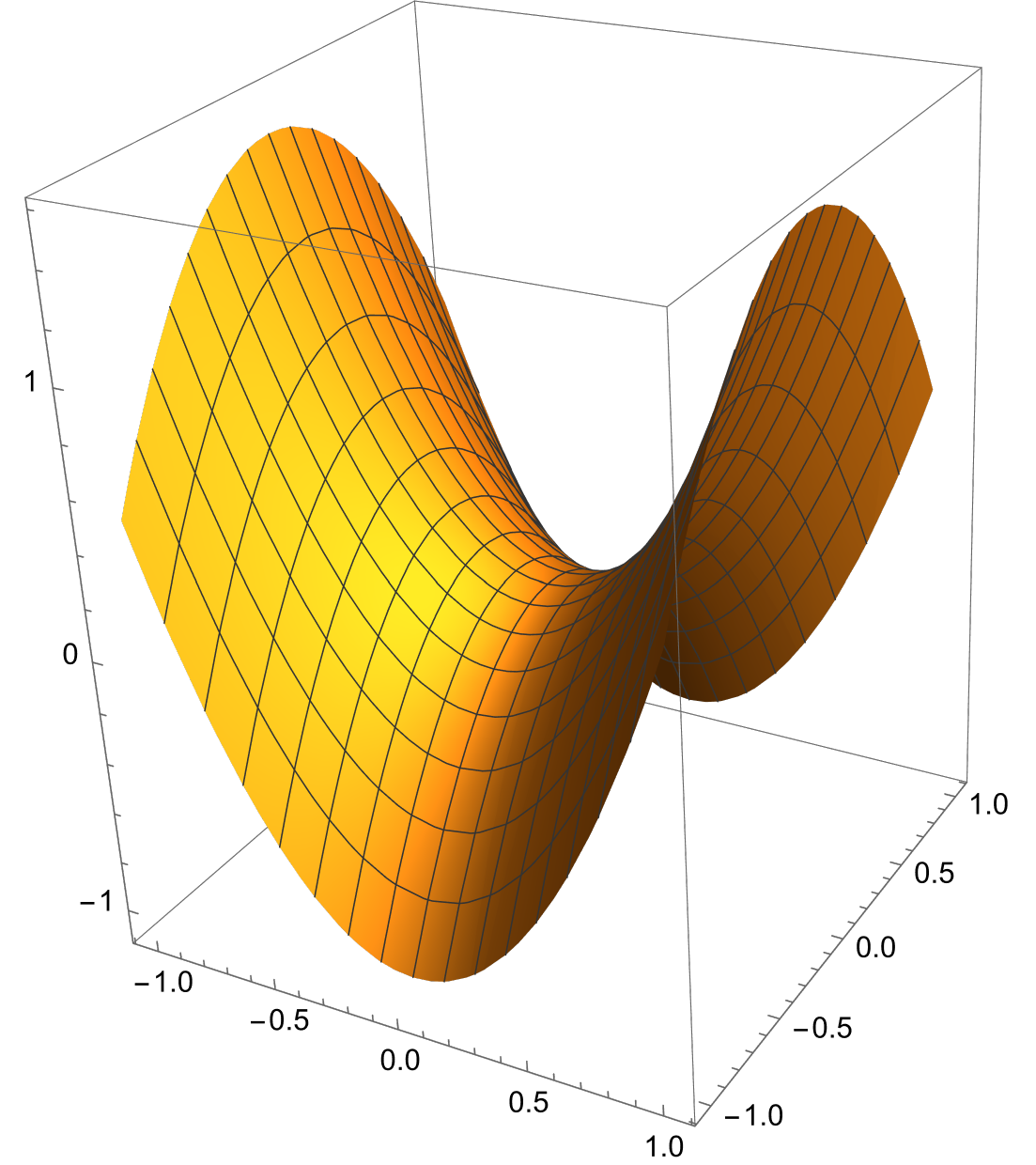}\quad \includegraphics[width=.31\textwidth]{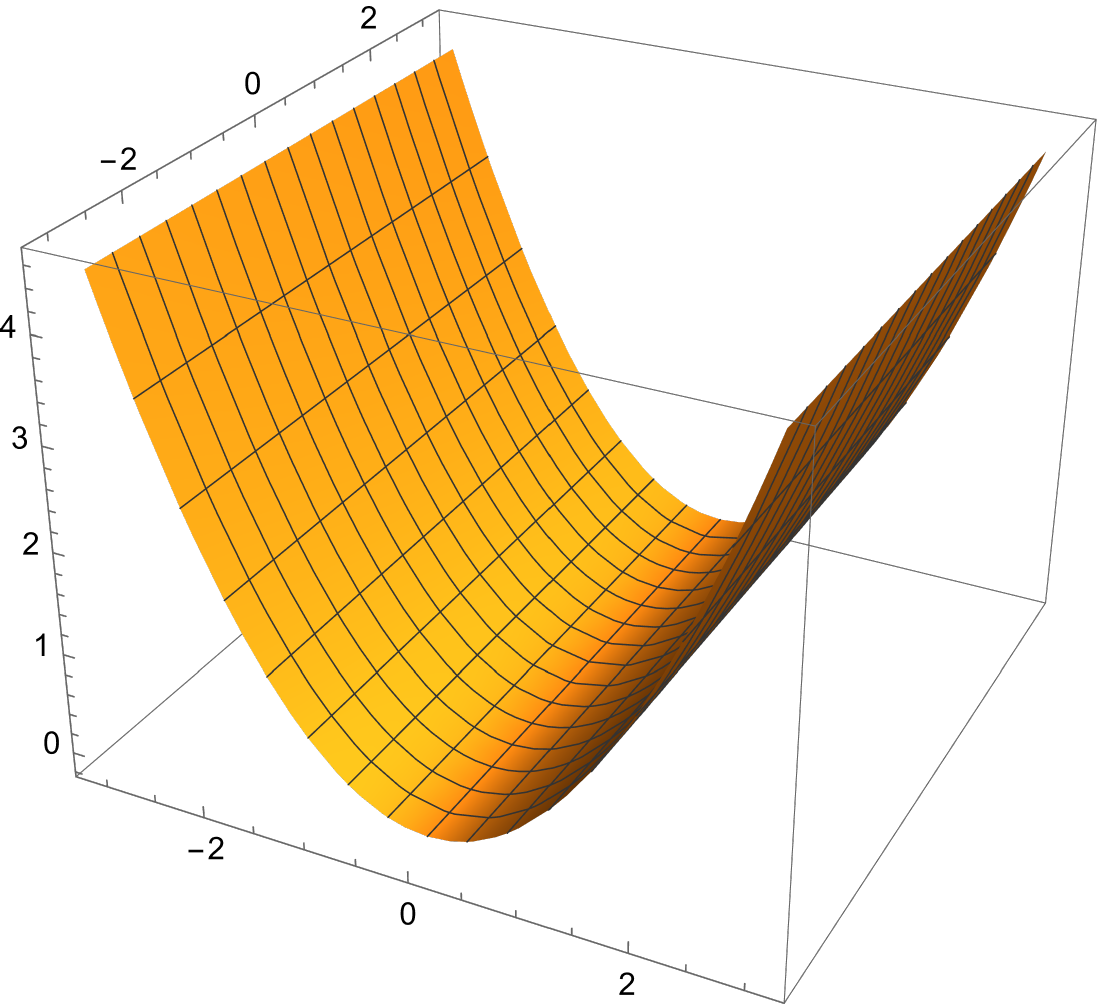}
\quad \includegraphics[width=.31\textwidth]{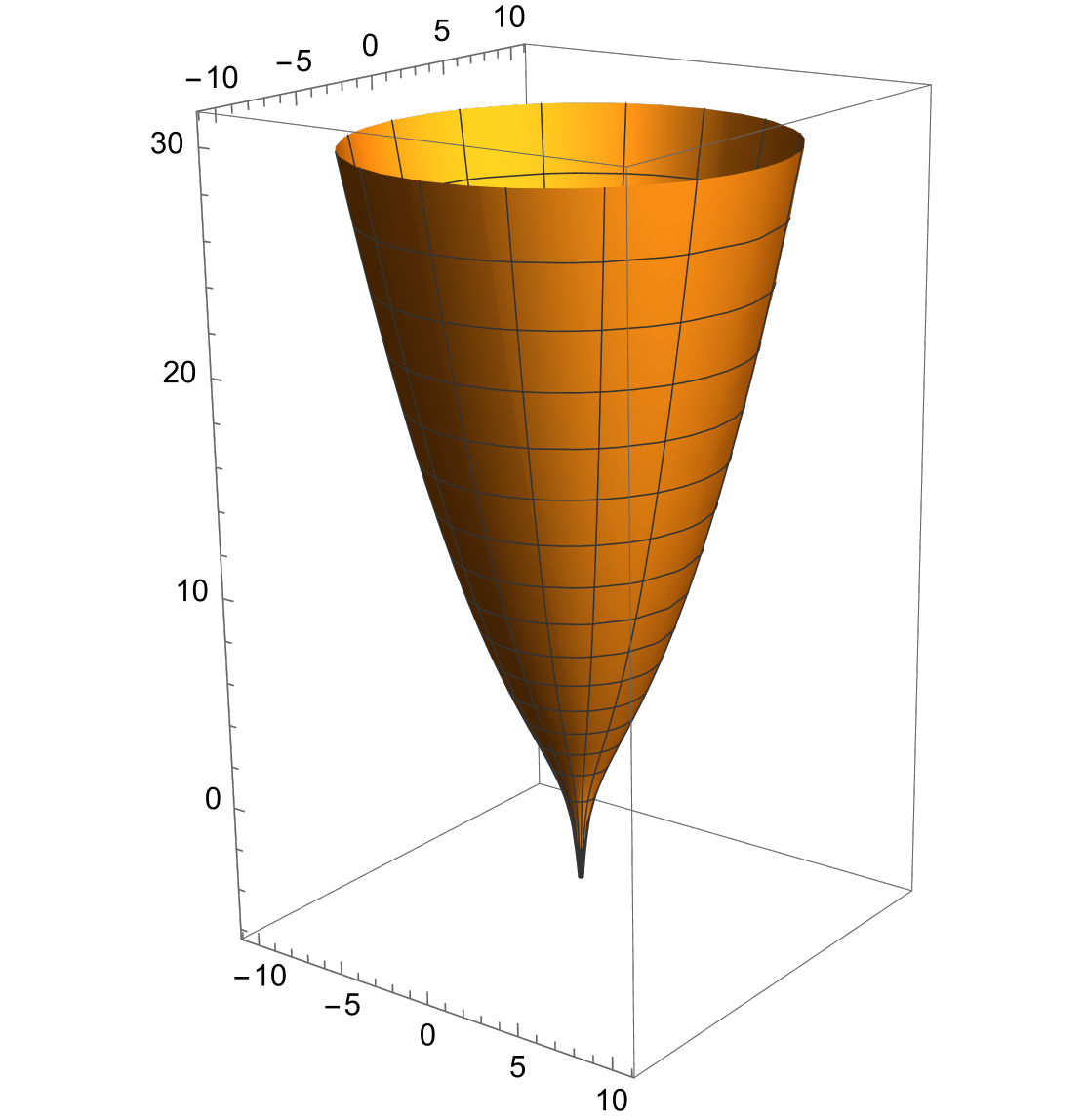} \end{center}
\caption{Separable surfaces of Prop. \ref{pr1}, case  $\Lambda\not=0$: translation surface (left), homothetical surface (middle), rotational surface. Here $\Lambda=1$. }\label{fig2}
\end{figure}

In this section, we prove Prop. \ref{pr1}, the proof of which is almost immediate. See Figs. \ref{fig1} and \ref{fig2}.
 
\begin{enumerate}
 
\item  If $z(x,y)=f(x)+g(y)$, then Eq. \eqref{eq1} becomes $f''+g''=\Lambda/2$. Since $f''$ depends only on $x$ and $g''$ on $y$, there exists $m\in\r$ such that 
$$f''-\frac{\Lambda}{2}=m=-g''$$
 and the solutions are given in \eqref{sol2}.   
\item If $z(x,y)=p(x)q(y)$, then Eq. \eqref{eq1} is $p''q+pq''=\Lambda/2$. Discarding the trivial cases $p=0$ or $q=0$ where the solution is a horizontal plane, we divide by $pq$, obtaining 
\begin{equation}\label{pq}
\frac{p''}{p}+\frac{q''}{q}=\frac{\Lambda}{2pq}.
\end{equation}
If $\Lambda=0$, then there exists $m\in\r$ such that 
$$\frac{p''}{p}=m=-\frac{q''}{q}.$$
If $m=0$, then $p$ and $q$ are linear functions. If $m=-k^2<0$, the solutions are given in \eqref{sol31}. If $m=k^2>0$, we obtain \eqref{sol31} again by interchanging the roles of $p$ and $q$.

Suppose $\Lambda\not=0$. Differentiating \eqref{pq} with respect to $x$ and then with respect to $y$, we obtain $ \Lambda p'q' =0$. Thus $p'=0$ or $q'=0$; that is, $p$ or $q$ are constant functions. If, for example,  $q(x)=q_0\not=0$ is constant, then \eqref{pq} becomes $p''q_0=\Lambda/2$ and the solution is given in \eqref{sol32}. 
\item Let $r=x^2+y^2$. Then  Eq. \eqref{eq1} becomes $4rh''+4h'=\Lambda/2$, whose solution is given in \eqref{sol4}.

\end{enumerate}
 
\begin{remark} If $\Lambda=0$, solutions \eqref{sol2} and \eqref{sol31} can be expressed in complex notation. If $w=x+iy\in\mathbb{C}$ and $\overline{w}$ denotes the conjugate of $w$, then the solutions   \eqref{sol2} and \eqref{sol31} can be written, respectively, as  
\begin{equation*}
\begin{split}
z(x,y)&=\mbox{Re}(\frac{m}{2}w^2+\alpha w+\beta),\\
 z(x,y)&=\alpha\cos w+\overline{\alpha}\cos\overline{w}+\beta\sin w+\overline{\beta}\sin\overline{w},
 \end{split}
 \end{equation*}
where $\alpha,\beta\in\mathbb{C}$.
\end{remark} 
%%%%%%%%%%%%%%%%%%%%%%%%%%%%%%%%
\section{First step of the proof of Theorem \ref{t1}}\label{s3}
%%%%%%%%%%%%%%%%%%%%%%%%

 We consider the general case in the implicit equation \eqref{s1}. From now on, we will discard the cases covered by Prop. \ref{pr1}.
 If $z=z(x,y)$, differentiating   with respect to $x$  we have 
$f'+h'z_x=0$ and differentiating again, we arrive at $f''+h''z_x^2+h' z_{xx}=0$. Doing the same with the variable $y$,   equation \eqref{eq1} becomes 
\begin{equation}\label{eq2}
(f''+g'') h'^2+(f'^2+g'^2)h''+\frac{\Lambda}{2} h'^{3}=0.
\end{equation}
We will obtain a first integral of the second-order differential equation \eqref{eq2}, yielding  three first-order differential equations \eqref{ss1} and \eqref{ss2} for the functions $f$, $g$ and $h$.

  Since the functions $f$, $g$ and $h$ are not constant, let us introduce   the following   variables 
\begin{equation*}
u=f(x),\quad v=g(y),\quad w=h(z).
\end{equation*}
 Then the implicit equation \eqref{s1} is simply 
 \begin{equation}\label{sep}
 u+v+w=0.
 \end{equation}
   Define the functions
$$X(u)=f'(x)^2,\quad Y(v)=g'(y)^2,\quad Z(w)=h'(z)^2.$$
In terms of these new functions, the cases discarded by Prop. \ref{pr1} are:
\begin{enumerate}
\item When one of the functions $X$, $Y$ or $Z$ is constant (cylindrical or  translation surfaces).
\item $Z(w)=e^{c_1w}$, with $c_1\not=0$ (homothetical).
\item $X(u)=a_1u+a_2u$ and $Y(v)=b_1v+b_2$ with $a_1=b_1$ (rotational).
\end{enumerate}
Differentiating the functions $X$, $Y$ and $Z$, we have
$$X'(u)= 2f''(x), \quad Y'(v) =2g''(y),\quad Z'(w) =2h''(z).$$
With this change of variables, equation (\ref{eq2}) takes the form 
\begin{equation}\label{eq0}
(X'+Y')Z+(X+Y)Z'+\Lambda Z^{3/2}=0,
\end{equation}
for all values of $u$, $v$ and $w$ under the condition \eqref{sep}.  

We now differentiate Eq. \eqref{eq0} with respect to $u$, $v$ and $w$. Notice that the variables $u$, $v$ and $w$ are not independent because  of \eqref{sep}. Thus, it is not possible, for example, to differentiate the left-hand side of \eqref{eq0} with respect to $u$ and ensure that it is $0$. However, we will use the following lemma, whose proof is immediate.

\begin{lemma} \label{le1}
Let $Q=Q(u,v,w)$ be a smooth function defined in a domain $\Omega\subset\r^3$. If $Q(u,v,w)=0$ for any triple of the section $\Omega\cap\Pi$, where $\Pi$ is the plane with equation $u+v+w=0$, then on this section we have $Q_u=Q_v=Q_w$, or equivalently,
\begin{equation}\label{l1}
Q_u-Q_v=0, \quad Q_u-Q_w=0,\quad Q_v-Q_w=0,
\end{equation}
where $Q_u$, $Q_v$ and $Q_v$ are the partial derivatives of $Q$ with respect to $u$, $v$ and $w$, respectively.
\end{lemma}

By applying \eqref{l1}, we differentiate \eqref{eq0} with respect to $u$ and $v$, obtaining 
\begin{equation}\label{eq3}
Z \left(X''-Y''\right)+Z' \left(X'-Y'\right)=0.
\end{equation}
 We apply \eqref{l1} again, differentiating  with respect to $u$ and $v$,
 \begin{equation}\label{eq32}
 Z \left(X^{(3)}+Y^{(3)}\right)+Z' \left(X''+Y''\right)=0.
 \end{equation}
Viewing    Eqs. \eqref{eq3} and \eqref{eq32} as a linear system in $Z$ and $Z'$ and since $Z\not=0$,   the determinant of the coefficients must vanish. This gives
 \begin{equation}\label{eq33}
 X''^2-Y''^2- (X'-Y' ) (X^{(3)}+Y^{(3)} ) =0.
 \end{equation}
In this equation, the variable $w$ does not appear. If we apply \eqref{l1} by differentiating \eqref{eq33} with respect to $u$ and $w$, then we  are  differentiating \eqref{eq33} with respect to $u$ and the   equation still  holds. The same applies to the variable $v$. Having established this,  we  only differentiate \eqref{eq33} with respect to $u$ and then with respect to $v$, obtaining
\begin{equation}\label{eq4}
X^{(4)}  Y'' -Y^{(4)}  X'' =0.
\end{equation}
Since the roles of $X$ and $Y$ are symmetric, the following arguments will be applied to $X$, and   similar reasoning holds  for the function $Y$. 

The remainder of this section is devoted to  the case  $X'' =0$ identically in \eqref{eq4}, while the  general case $X''\not=0$ will be considered in Sect. \ref{s4}.   

Assume  $X''=0$. Then $X(u)=a_1u+a_2$, where  $a_1\not=0$. Substituting this  into \eqref{eq33}, we have 
\begin{equation}\label{y2}
(Y'-a_1)Y'''-Y''^2=0.
\end{equation}
In order to solve this equation, we separate the cases where $Y'-a_1=0$ identically and $Y'-a_1\not=0$.
\begin{enumerate}
\item Case  $Y'-a_1=0$ identically. The solution to this equation is  $Y(v)=a_1v+b_2$. But this implies that the surface is rotational because $X$ and $Y$ are linear functions of their variables and the coefficient of $v$ coincides with that of $u$ in the expression for $X(u)$. This case was initially discarded.

\item Case  $Y'-a_1\not=0$. The solution to \eqref{y2}  is 
$$Y(v)=a_1 v+\frac{1}{b_1}e^{b_1 v+b_2}+b_3,\quad b_1\not=0.$$
Now \eqref{eq32} becomes 
 $$b_1 e^{b_1 v+b_2}  (Z'+b_1Z)=0.$$
 Thus $Z'+b_1Z=0$, whose solution is $Z(w)=c_1e^{-b_1 w}$, with$c_1\not=0$. However, this situation corresponds with the homothetical surfaces in Prop. \ref{pr1}, which were also discarded from the beginning. In fact, when $Z=e^{c_1 w}$  in Prop. \ref{pr1}, the function $Y$ is constant, which is not possible here because $Y'-a_1\not=0$. Continuing with our current argument,  and using the value of $Z(w)$,  identity \eqref{eq0} becomes 
 $$c_1e^{-b_1w}\left(a_2 b_1+b_1 b_3-a_1(2+b_1w+\Lambda\sqrt{c_1}e^{-b_1w/2})\right)=0.$$
  Since the functions $\{1,w,e^{-b_1w/2}\}$ are linearly independent, we deduce that $ a_1b_1  =0$,  which is a contradiction.
\end{enumerate}

%%%%%%%%%%%%%%%%%%%%%%%%%%%%%%%%
\section{Second step of the proof of Theorem \ref{t1}: case $X''\not=0$}\label{s4}
%%%%%%%%%%%%%%%%%%%%%%%%

In this section, we complete the proof of Thm. \ref{t1}. Assume  $X''\not=0$ and, by symmetry,  we also assume $Y''\not=0$.  From \eqref{eq4}, there exists a real constant $m\in\r$ such that 
  \begin{equation}\label{mm}
  \frac{X^{(4)}}{X''}=m=\frac{Y^{(4)}}{Y''}.
  \end{equation}
  The discussion depends on the sign of $m$.
 
   \subsection{Case  $m=0$.} 

We prove that this case is not possible. The solutions to \eqref{mm} are
\begin{equation*}
\begin{split}
X(u)&=a_1u^3+a_2 u^2+a_3 u+a_4,\\
Y(v)&=b_1 v^3+b_2 v^2+b_3 v+b_4.
\end{split}
\end{equation*}
 Moreover, $a_1^2+a_2^2\not=0$ and $b_1^2+b_2^2\not=0$ to ensure $X''\not=0$ and $Y''\not=0$. Substituting into \eqref{eq33}, we obtain a polynomial equation of the type 
  $A_0+A_1 u+B_1v+A_2 u^2+B_2 v^2=0$, where
  \begin{equation*}
  \begin{split}
   A_2&=18 a_1(a_1-b_1),\\
    B_2&=18 b_1(a_1-b_1),\\
     A_1&=12 a_2(a_1-b_1),\\
      B_1&=12b_2(a_1-b_1),\\
       A_0&=4(a_2^2-b_2^2)-6(a_1+b_1)(a_3-b_3-b_4).
       \end{split}
       \end{equation*}
From $A_2=0$, we have two cases.
\begin{enumerate}
\item Case $a_1=0$. Then $B_2=0$ implies $b_1=0$ and from $A_0=0$, we have $b_2^2=a_2^2$ with $a_2\not=0$. Thus $b_2=\pm a_2$ and now  \eqref{eq32} becomes   $\pm4a_2Z'=0$, which is a contradiction because $a_2\not=0$ and $Z$ is not a constant function. 
\item Case $a_1=b_1$ and $a_1\not=0$.  Then $A_0=0$ is written as   
  $$a_2^2-b_2^2+3a_1(b_4+b_3-a_3)=0.$$
  Now \eqref{eq32} becomes 
  $$6 a_1 Z+(a_2+b_2+3a_1(u+v))Z'=0.$$
  Since $u+v=-w$,   this equation is converted to
\begin{equation*} 
6 a_1 Z+(a_2+b_2-3a_1 w)Z'=0.
\end{equation*}
Given $a_1\not=0$, the solution  is
$$Z(w)=c_1(a_2+b_2-3a_1 w)^2,\quad c_1\not=0.$$
  Substituting this into \eqref{eq0} we obtain a polynomial equation in the variables $u$, $v$ and $w$. When we replace $w=-u-v$, this equation can be expressed as 
  $$P_1+\Lambda(c_1(a_2+b_2+3a_1 (u+v))^{3/2}=0.$$
  This implies
  $$P_1^2-\Lambda^2(c_1(a_2+b_2+3a_1 (u+v))^{3 }=0.$$
  This identity is a   polynomial equation $P(u,v,u^rv^s)=0$ with $r+s\leq 8$.   Thus, all coefficients must vanish.  The coefficient of $u^8$ is $81a_1^6 c_1^2$. However, $81a_1^6 c_1^2=0$ is a contradiction because $a_1,c_1\not=0$.

  \end{enumerate}

%%%%%%%%
  \subsection{Case  $m=-k^2<0$.} 
  %%%%%%%%
  
   For this value of $m$, the solutions  for $X$ and $Y$ in \eqref{mm} are
  $$X(u)=a_1\cos (ku)+a_2\sin(ku)+a_3 u+a_4,$$
   $$Y(v)=b_1\cos (kv)+b_2\sin(kv)+b_3 v+b_4.$$
   Here $a_1^2+a_2^2\not=0$ and $b_1^2+b_2^2\not=0$ in order to ensure again that $X''\not=0$ and $Y''\not=0$. 
   Substituting these into \eqref{eq33}, we obtain
   $$ k \left(a_1^2+a_2^2-b_1^2-b_2^2\right)+(b_3-a_3) (a_1 \sin (k u)-a_2 \cos (k u)+b_1 \sin (k v)-b_2 \cos (k v))=0.$$
 Since the trigonometric functions are linearly independent and $k\not=0$, we deduce
 \begin{equation*}
 \begin{split}
 b_2(a_3-b_3)&=0,\\
  a_2(a_3-b_3)&=0,\\
   b_1(a_3-b_3)&=0,\\
 a_1(a_3-b_3)&=0,\\
 a_1^2+a_2^2-b_1^2-b_2^2&=0.
 \end{split}
 \end{equation*}
 It is immediate that $a_3-b_3=0$ because otherwise, $a_1=b_1=a_2=b_2=0$, which is not possible. Letting $b_3=a_3$, the last equation implies $a_1^2+a_2^2=b_1^2+b_2^2$. Since this number is non-zero,  there exist $r\not=0$ and $d_1,d_2$ such that 
 \begin{equation*}
 \begin{split}
 &a_1=r\cos d_1,\quad a_2=r\sin d_1,\\
 &b_1=r\cos d_2,\quad b_2=r\sin d_2.
 \end{split}
 \end{equation*}
 
 Let $v=-u-w$. By computing identity \eqref{eq32}, we have
$$ \sin (\frac{1}{2} (k (2 u+w)-d_1+d_2)) \left(Z' \cos (\frac{1}{2} (k w+d_1+d_2))+k Z \sin (\frac{1}{2} (k w+d_1+d_2))\right)=0.$$ 
Since $\sin (\frac{1}{2} (k (2 u+w)-d_1+d_2))\not=0$ because $u$ and $w$ are here independent variables, we deduce that   the term in parentheses in the above equation must vanish identically. This leads to 
$$\tan (\frac{1}{2} (k w+d_1+d_2))=-\frac{Z'}{kZ}.$$
The solution to this equation is  $$Z(w)=c_1(1+\cos( kw+d_1+d_2 ) ),\quad c_1\not=0.$$
 Now \eqref{eq3} is trivially $0$ and \eqref{eq0} becomes $P_1+\sqrt{2}\Lambda P_2=0$, where
 \begin{equation*}
 \begin{split}
 P_1&=a_3c_1\left(2 +2\cos(kw+d_1+d_2)+kw\,\sin(kw+d_1+d_2)\right),\\
 P_2&= (c_1(1+\cos(kw+d_1+d_2))^{3/2}.
 \end{split}
 \end{equation*}
 Thus, 
 \begin{equation}\label{p1}
 P_1^2- 2\Lambda^2c_1^3 (1+\cos(kw+d_1+d_2))^3=0.
 \end{equation}
  This equation is a polynomial equation on the functions  $w^{n_1} \sin(n_2(kw+d_1+d_2) )  \cos(n_3(kw+d_1+d_2))  $, where $n_i\in\mathbb{N}$, $n_1\leq 2$ and $n_2,n_3\leq 3$. Thus, all coefficients vanish identically. The coefficient of $\cos (3(kw+d_1+d_2))$ in \eqref{p1} is $-\frac12c_1^3\Lambda^2$. Since $c_1\not=0$, we have $\Lambda=0$. As a consequence, the case $\Lambda\not=0$ cannot occur. This proves, together with the same conclusion in the case $m=k^2>0$ below, the situation $\Lambda\not=0$ described in Thm. \ref{t1}, where only the separable surfaces of Prop. \ref{pr1} are possible. 
  
  We now continue the reasoning assuming from now on that $\Lambda=0$. The computation of the coefficient of $w\sin (2(kw+d_1+d_2))$  in \eqref{p1} is $2a_3^2c_1^2k$. Thus, $a_3=0$ and, consequently, $b_3=0$. Now \eqref{eq0} reduces to 
  $$4k^2c_1^2(a_4+b_4)^2(\sin(2(kw+d_1+d_2)))^2=0.$$
  This implies $b_4=-a_4$. Returning to the previous expressions, we have
  \begin{equation*}
  \left\{\begin{split}
  X(u)&=r\cos (ku-d_1) +a_4,\\
  Y(v)&=r\cos (kv-d_2)-a_4,\\
   Z(w)&=c_1 (1+\cos( kw+d_1+d_2)).
\end{split}
\right.
\end{equation*}
Substituting these into \eqref{eq0} results in a trivial identity. Thus, these are the solutions that we seek, which coincide with \eqref{ss1}.

%%%%%%%%
  \subsection{Case  $m=k^2>0$.}  
  
  The arguments   are similar to the case $m<0$. We deduce     $\Lambda=0$ and  
   \begin{equation*}
    \left\{
    \begin{split}
  X(u)&=r\cosh (ku-d_1) +a_4,\\
  Y(v)&=r\cosh (kv-d_2)-a_4,\\
   Z(w)&=c_1 (1+\cosh( kw+d_1+d_2)).
\end{split}
\right.
\end{equation*}
These solutions coincide with \eqref{ss2}.
 
  %%%%%%%%
  \section{Explicit examples of separable surfaces}\label{s5}
  %%%%%%%%
  
  In this section, we provide some explicit examples of separable surfaces  when $\Lambda=0$ by solving the  systems     \eqref{ss1}  and \eqref{ss2} in Thm. \ref{t1}. Note that  numerous examples of surfaces can be obtained through various   choices of  the integration constants in both systems. To be precise, we present three examples. In the first two examples,  the function $h(z)$ is explicitly obtained. When solving the systems \eqref{ss1} and \eqref{ss2} for the functions $f$ and $g$, elliptic integrals of the first kind arise, and the   solutions can be expressed using the inverse of these functions.    In the third example, we solve explicitly \eqref{ss2} for suitable choices of constants, obtaining the parametrization of the surface. 
  
We derive  the symmetry properties of these surfaces using the Schwarz reflection principle for harmonic functions. Since $\Lambda=0$, the surface is the graph of a function $\varphi=\varphi(x,y)$ which is harmonic on the $xy$-plane. Consequently, we   have the following two results, which are analogous to those appearing   in the theory of minimal surfaces \cite{hk}. 
\begin{enumerate}
\item If $\Sigma$ contains a horizontal line $L$, then $\Sigma$ can be extended analytically by     $180^0$ rotations across $L$. Equivalently, if $\Sigma$ satisfies Eq. \eqref{eq1} for $\Lambda=0$ and $L$ is a straight line contained in $\Sigma$, then $\Sigma$ is symmetric with respect to a $180^0$ rotation about $L$.
\item If $\Sigma$ meets   a vertical plane $\Pi$ orthogonally along a regular curve $C\subset\Sigma\cap\Pi$,   then $\Sigma$ can be extended analytically by symmetry  across $\Pi$. Equivalently, if $\Sigma$ satisfies Eq. \eqref{eq1}  for $\Lambda=0$  and $\Pi$ is a plane that meets  $\Sigma$ orthogonally, then  $\Sigma$ is symmetric under  reflection  across $\Pi$.
\end{enumerate}

%The difference with the minimal surfaces is that now   the line $L$ must be horizontal and that the plane $\Pi$ must be vertical. This ensures that the reflections preserve the harmonicity of the third coordinate $z(x,y)$ of the surface. 

%%%%%%
\subsection{First example}
%%%%%%%%

We begin with the system  \eqref{ss1}. Let $k=-2$.   Take $r=c_1=1$ and $d_1=d_2=a=0$. Let $\Sigma_1$ be the corresponding separable surface defined by the implicit equation \eqref{s1}. Then we have
     \begin{equation*}
     \left\{
  \begin{split}
 f'(x)^2&=\cos   (2f(x)),\\
g'(y)^2&=\cos  (2g(y)),\\
h'(z)^2&=1+\cos (2h(z)).
\end{split}
\right.
\end{equation*} 
 The first equation yields
   $$\int\frac{df}{\sqrt{\cos (2 f)}}=\int\frac{df}{\sqrt{1-2(\sin f)^2}}= x+\lambda,$$
  with $\lambda\in\r$.  For simplicity, let $\lambda =0$. This integral is the elliptic integral of the first kind $F=F(\phi,2)$. Thus,  $F(f,2)=x$ and, analogously, $F(g,2)=y$.  Let $M=M(\phi,2)$ be a branch of the inverse function of $F$ with $M(0,2)=0$.   Thus, 
 \begin{equation*}
 \left\{ \begin{split}
 f(x)&=M(x,2),\\
 g(y)&=M(y,2).
 \end{split}\right.
 \end{equation*}
 
\begin{figure}[hbtp]
\begin{center}\includegraphics[width=.45\textwidth]{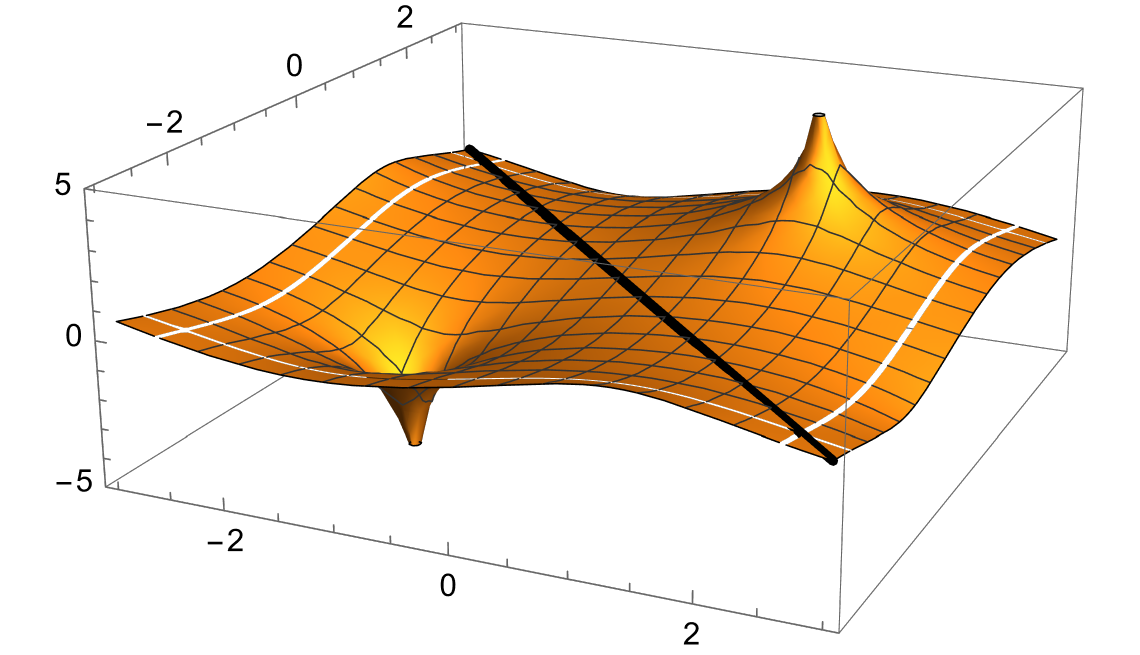}\quad \includegraphics[width=.4\textwidth]{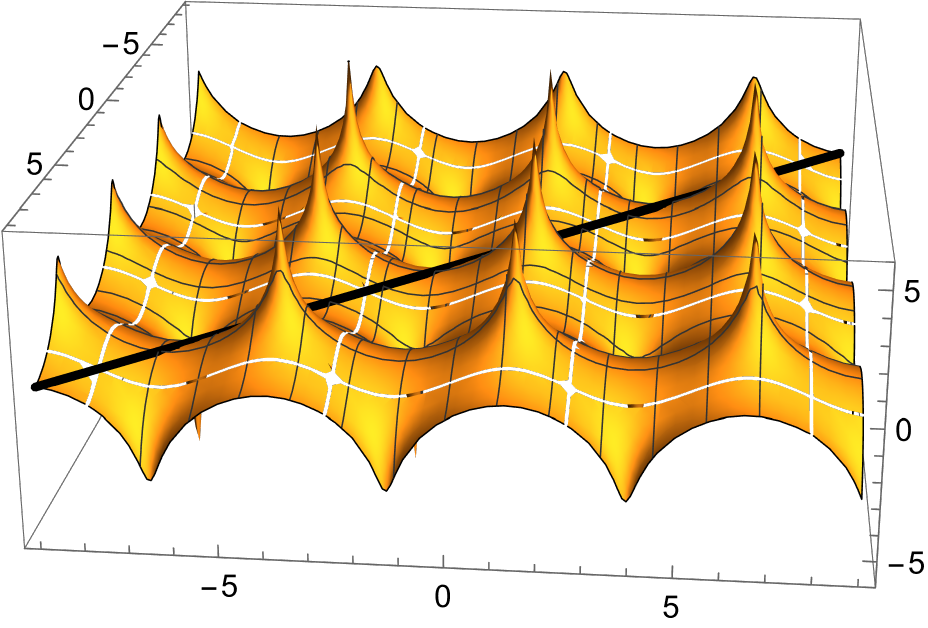}
\ \end{center}
\caption{The surface $\Sigma_1$ defined in \eqref{ex1}. Left: a piece of the surface showing a straight line (black) contained in the surface. Right: the same surface after reflection of $180^0$ degrees across these lines.   }\label{fig3}
\end{figure}

 We now find the function $h=h(z)$. In this case,  we   find $h$ because 
$$\int\frac{dh}{\sqrt{1+\cos(2h)}}=\frac{1}{\sqrt{2}} \tanh ^{-1}(\sin (h)) .$$
Thus, up to an additive constant, 
$$h(z)=\sin^{-1}(\tanh(\sqrt{2}z)).$$
The surface is then defined by   the implicit equation:
\begin{equation}\label{ex1}
\Sigma_1: M(x,2)+M(y,2)+\sin^{-1}(\tanh(\sqrt{2}z))=0.
\end{equation}
A picture of the surface is shown in   Fig. \ref{fig3}, left. The function $f(x)=M(x,2)$ (also $g=g(y)$) is an odd function.   Since $M(-x,2)=-M(x,2)$, the straight line $L=\{x,-x,0):x\in\r\}$ is contained in $\Sigma$.  Therefore   $\Sigma_1$    is invariant under a    $180^0$  rotation about $L$. This allows the surface to be  extended by successive reflections. In fact, the surface contains two families of parallel straight lines $\{L_m,L^m:m\in\mathbb{Z}\}$, all of which are contained in the horizontal plane $z=0$ and  defined by:
$$L_m:=\{(x,-x+4m,0):x\in\r\},\quad L^m:=\{(x,x+4m-2,0):x\in\r\}.$$

\begin{figure}[hbtp]
\begin{center}\includegraphics[width=.4\textwidth]{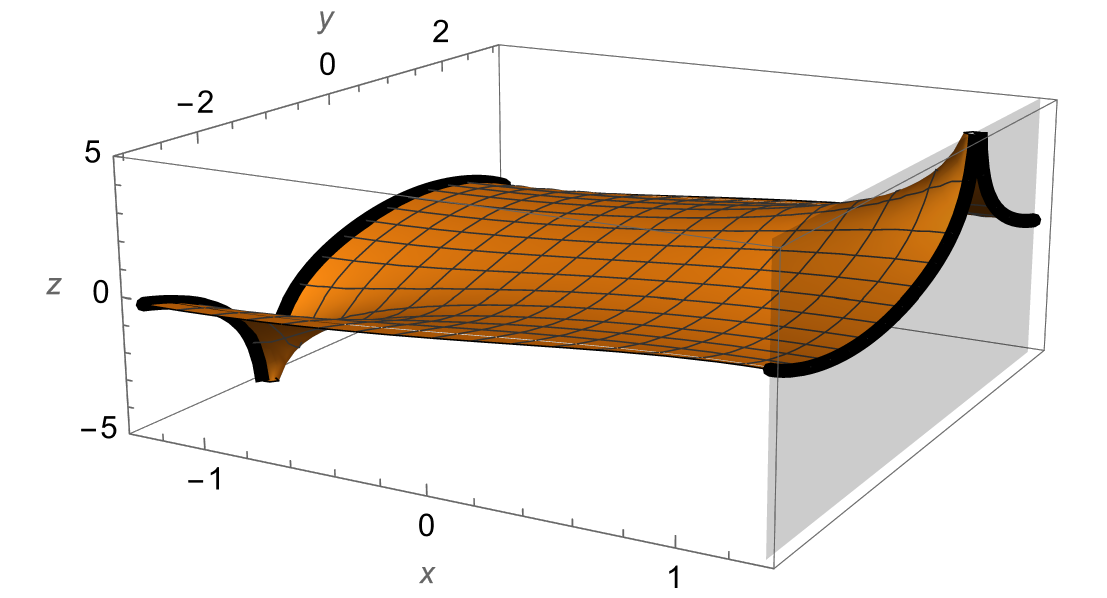}\quad \includegraphics[width=.5\textwidth]{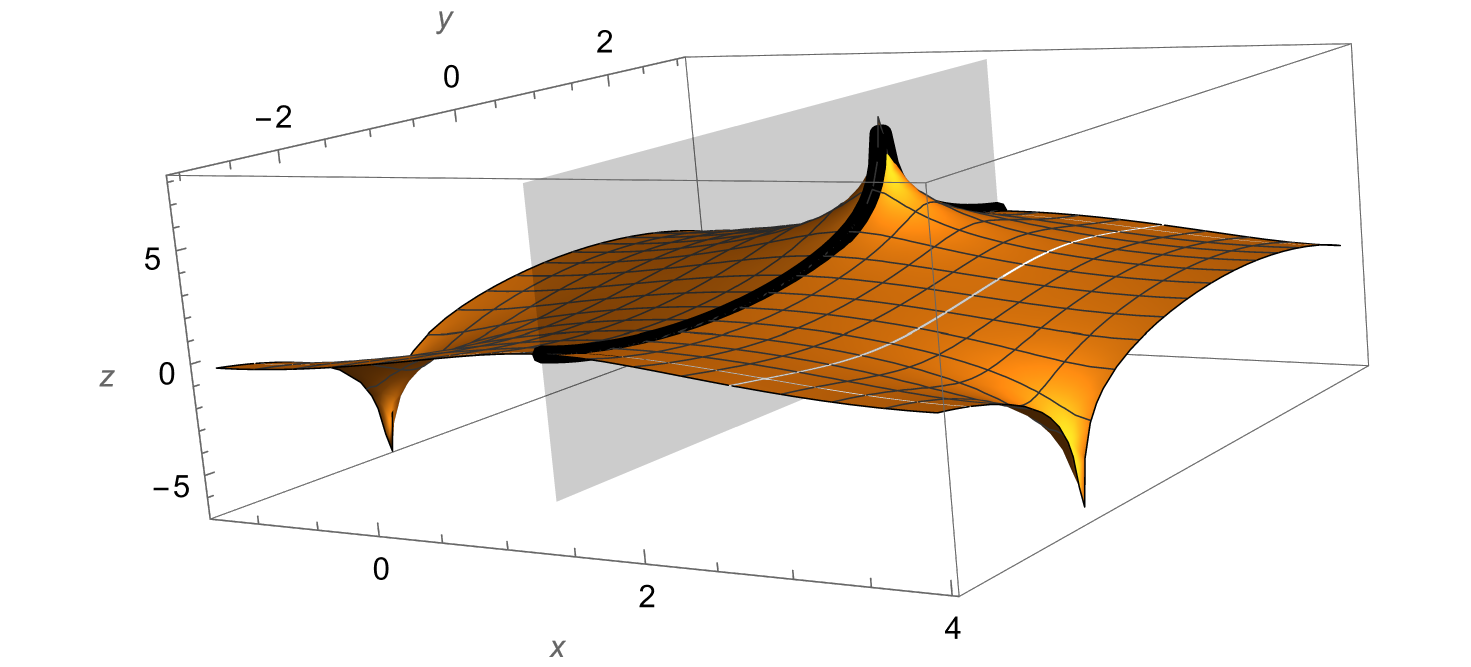}
\ \end{center}
\caption{The surface $\Sigma_1$ defined in \eqref{ex1}. Left: a piece of the surface showing the intersection curves (black) of $\Sigma_1$ with the vertical planes of equations $x=\pm x_0$.   Right: the same surface after reflection across the plane $x=x_0$.     }\label{fig4}
\end{figure}

Other symmetries of the surface can be obtained by reflections across vertical planes.   In the case of the surface \eqref{ex1}, these planes occur at the curves of $\Sigma$ where $f'(x)=0$ and $g'(y)=0$. The first point where the derivative of $f$ vanishes is $x_0 \approx 1.311$. By the symmetry of $f$, the same occurs at $-x_0$.   Let $C=\Sigma_1\cap \{(x,y,z)\in\r^3: x=\pm x_0\}$. Since the normal $\nu$ of $\Sigma_1$ at $(x,y,z)\in\Sigma_1$ is proportional to $(f'(x),g'(y),h'(z))$, then $\nu$ is contained in the vertical planes defined by $x=\pm x_0$. This proves that $\Sigma_1$ meets both planes orthogonally (Fig. \ref{fig4}, left). Thus the reflections across these planes extend the surface $\Sigma_1$ (Fig. \ref{fig4}, right). Similarly, since the function $g$ coincides with $f$, the surface $\Sigma_1$ meets  the planes $y=\pm x_0$ orthogonally and, thus, $\Sigma_1$ can be also extended by reflections across these planes (Fig. \ref{fig4}, right). 

The above process of reflection across vertical planes can be repeated. These planes are situated at the points where $z$ approaches $\pm\infty$. The points $(x,y)$ where $z$ is $\infty$ are $\{((4m+1)x_0,(4n+1)x_0):m,n\in\mathbb{Z}\}$ and,  analogously, the points where $z$ is $-\infty$ are $\{((4m-1)x_0,(4n-1)x_0):m,n\in\mathbb{Z}\}$. Thus, the vertical planes of symmetry are given by $\{x=(2n-1)x_0\}$ and $\{y=(2n-1)x_0\}$. As a final observation, we point out that the surface is a global graph on $\r^2$ except for a discrete set points which is invariant under translations of $\r^2$ in two linearly independent directions.

\begin{proposition} If $\Sigma_1$ is the surface defined in \eqref{ex1}, then $\Sigma_1$ is doubly periodic. Moreover, $\Sigma_1$ is invariant under reflections across the lines $L_m$ and $L^m$ and under symmetries with respect to the planes  $\{x=(2m-1)x_0\}$ and $\{y=(2m-1)x_0\}$, $m\in\mathbb{Z}$. The surface is a global graph on $\r^2$ except for a discrete set points.

\end{proposition}

%%%%%%
\subsection{Second example}
%%%%%%%%
For the second example, we take $k=2$ in \eqref{ss2}. Again, let $r=c_1=1$ and $d_1=d_2=a=0$. The system to solve is     \begin{equation*}
\left\{  \begin{split}
 f'(x)^2&=\cosh   (2f(x)),\\
g'(y)^2&=\cosh  (2g(y)),\\
h'(z)^2&=1+\cosh (2h(z)).
\end{split}
\right.
\end{equation*} 
Again, the function $h$ can be explicitly found, obtaining
$$h(z)=\sinh^{-1}(\tan(\sqrt{2}z)).$$
For the function $f$ (and similarly for $g$), we have 
$$\int\frac{df}{\sqrt{1+2(\sinh f)2}}=x+\lambda.$$
Let $\lambda=0$. In this case,  $-i\, F(if,2)=x$ and $-i\, F(ig,2)=y$, and thus:
 \begin{equation*}
 \left\{ \begin{split}
 f(x)&=i\, M(ix,2),\\
 g(y)&=i\, M(iy,2).
 \end{split}\right.
 \end{equation*}

Let  $\Sigma_2$ be the corresponding surface, 
\begin{equation}\label{ex2}
\Sigma_2: f(x)+g(y)+ \sinh^{-1}(\tan(\sqrt{2}z))=0.
\end{equation}
The function $f=f(x)$ is an odd function defined on the interval $(-x_1,x_1)$, where $x_1\approx 1.31$ with the property that $\lim_{x\to\pm x_1}f(x)=\mp \infty$. Notice that $x_1$ is the same value $x_0$ as the function $f(x)=M(x,2)$ for $\Sigma_1$. Let $y_1=x_1$ and $z_1=\frac{\pi}{2\sqrt{2}}$ such that  $(x_1,y_1,z_1)\in\Sigma_2$

\begin{figure}[hbtp]
\begin{center}\includegraphics[width=.42\textwidth]{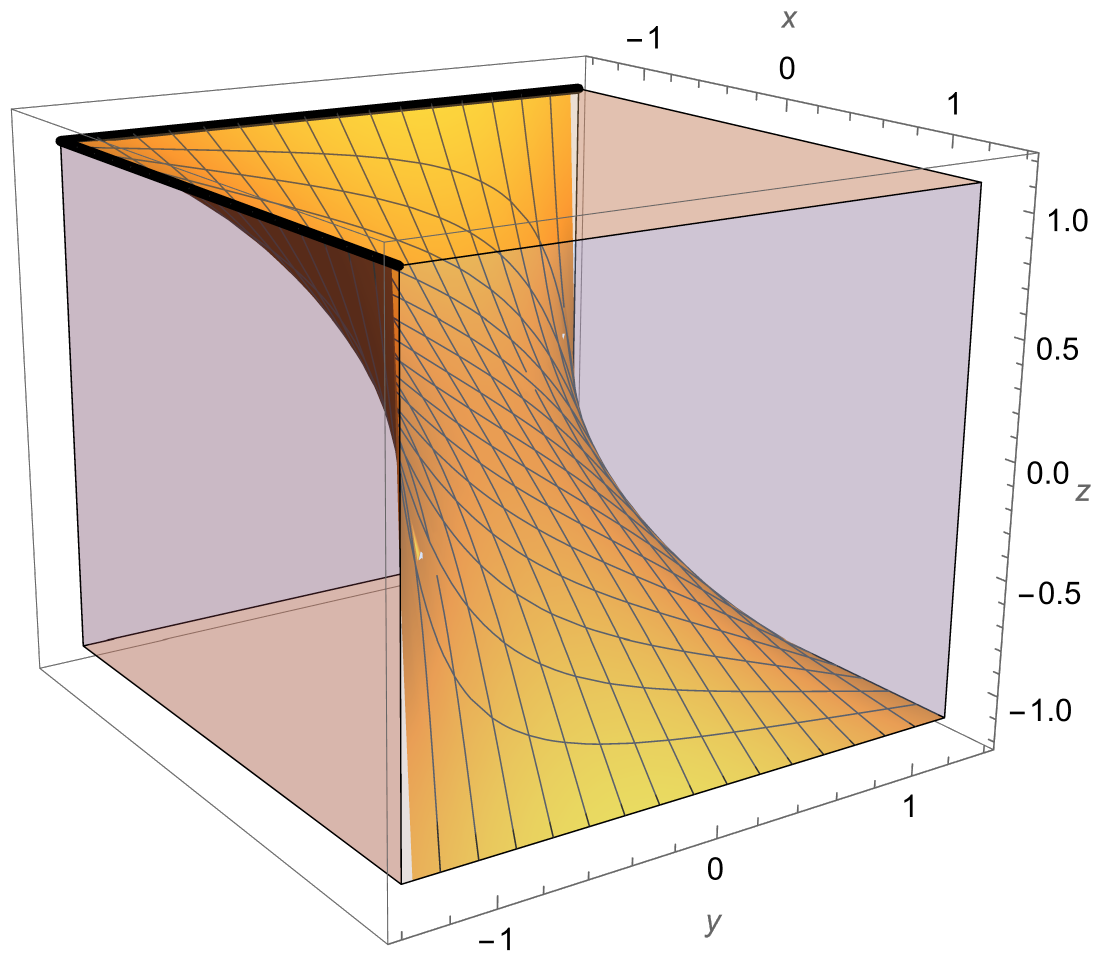}\qquad \includegraphics[width=.4\textwidth]{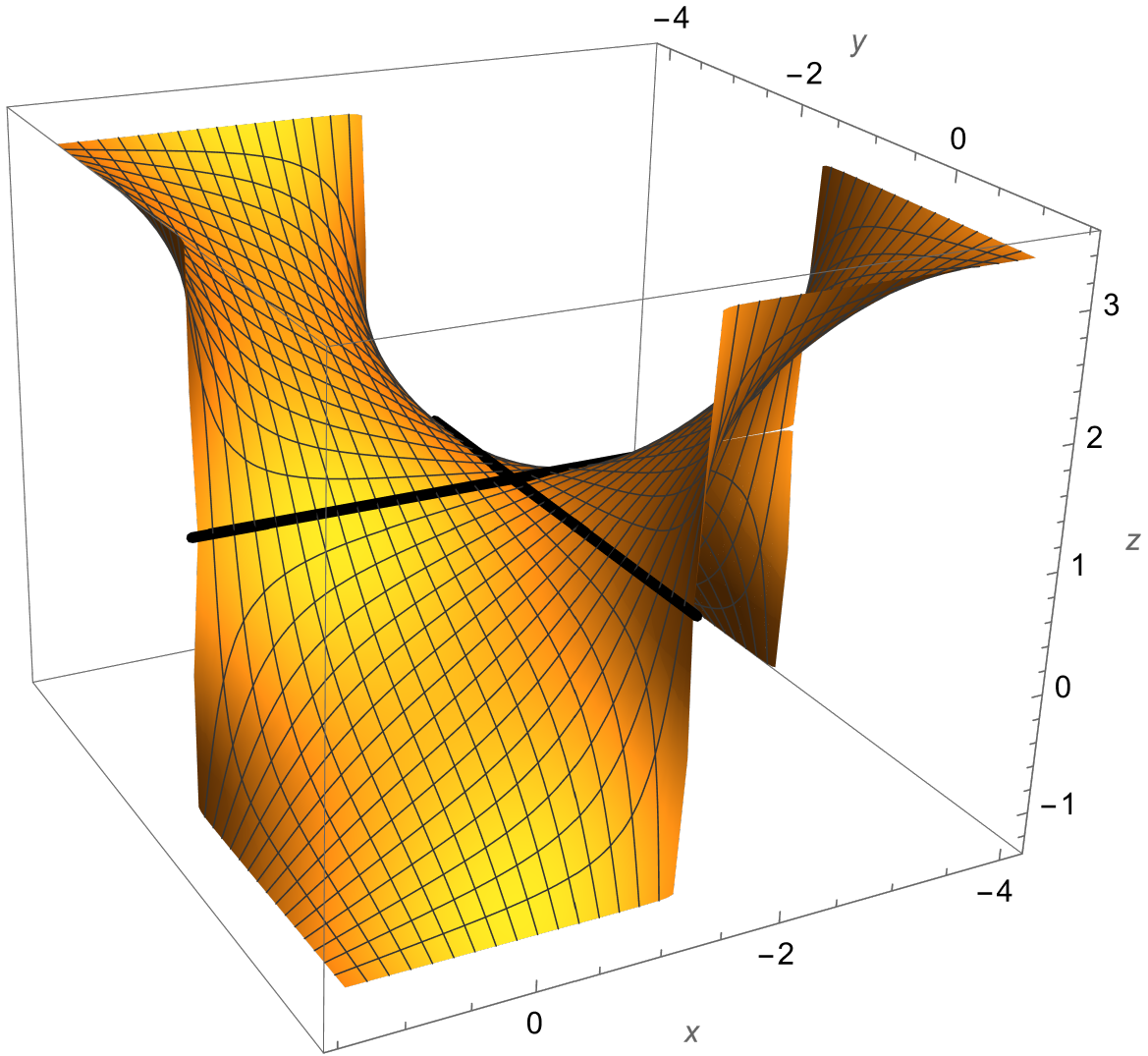}
 \end{center}
\caption{ The surface $\Sigma_2$  defined in  \eqref{ex2}. Left: a piece of the surface showing the intersection (black) of $\Sigma_2$ with the cuboid and the plane $z=z_1$.    Right: the same surface after successive reflections across horizontal lines.  }\label{fig5}
\end{figure}

The surface $\Sigma_2$ meets the cuboid with vertices $(\pm x_1,\pm y_1,\pm z_1)$ along six segments (Fig. \ref{fig5}, left). Only four of them are horizontal, namely, 
\begin{equation*}
\begin{split}
L_1&=\{(-x_1,y,z_1:-y_1<y<y_1\},\\
L_2&=\{(x,-y_1,z_1:-x_1<x<x_1\},\\
L_3&=\{(x_1,y,-z_1:-y_1<y<y_1\},\\
L_4&=\{(x,y_1,-z_1:-x_1<x<x_1\}.\\
\end{split}
\end{equation*}
Lines $L_1$ and $L_2$ are depicted in Fig. \ref{fig5}, left, as black lines. Thus, using the Schwarz reflection principle, the surface $\Sigma_2$ can be extended by a rotation of $180^0$ across $L_i$, $1\leq i\leq 4$ (Fig. \ref{fig5}, right). Since  $\Sigma$ can be reflected along two directions in space, then   $\Sigma_2$ is a doubly   periodic surface. Notice that although $\Sigma$ contains two vertical lines, the reflection across that lines does not imply that the surface can be extended along that direction. Consequently, the surface is a graph on $\r^2$ except for a discrete set of  points. 

\begin{proposition} If $\Sigma_2$ is the surface  defined in  \eqref{ex2}, then $\Sigma_2$ is doubly periodic.  The surface is a global graph on $\r^2$ except for a discrete set points.
\end{proposition}

%%%%%%
\subsection{Third example}
%%%%%%%%

\begin{figure}[hbtp]
\begin{center}\includegraphics[width=.5\textwidth]{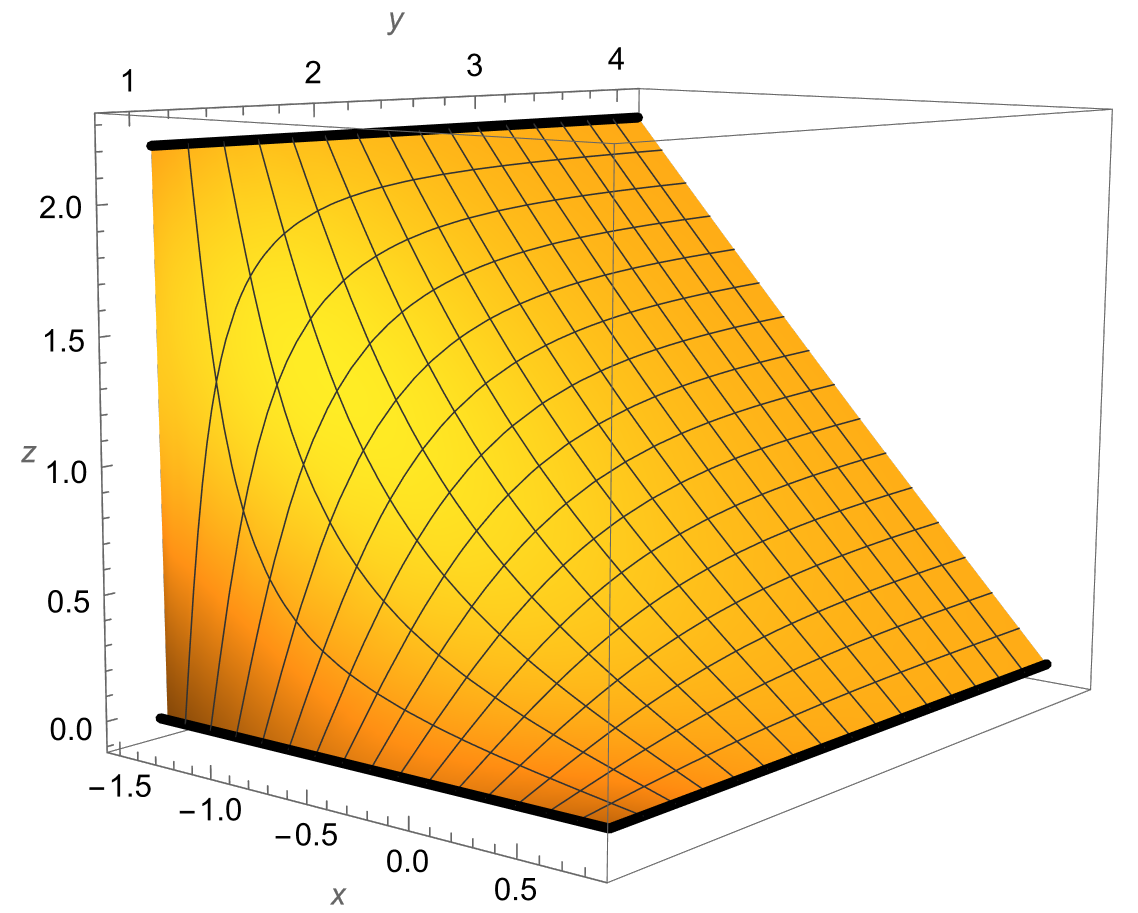}\qquad \includegraphics[width=.38\textwidth]{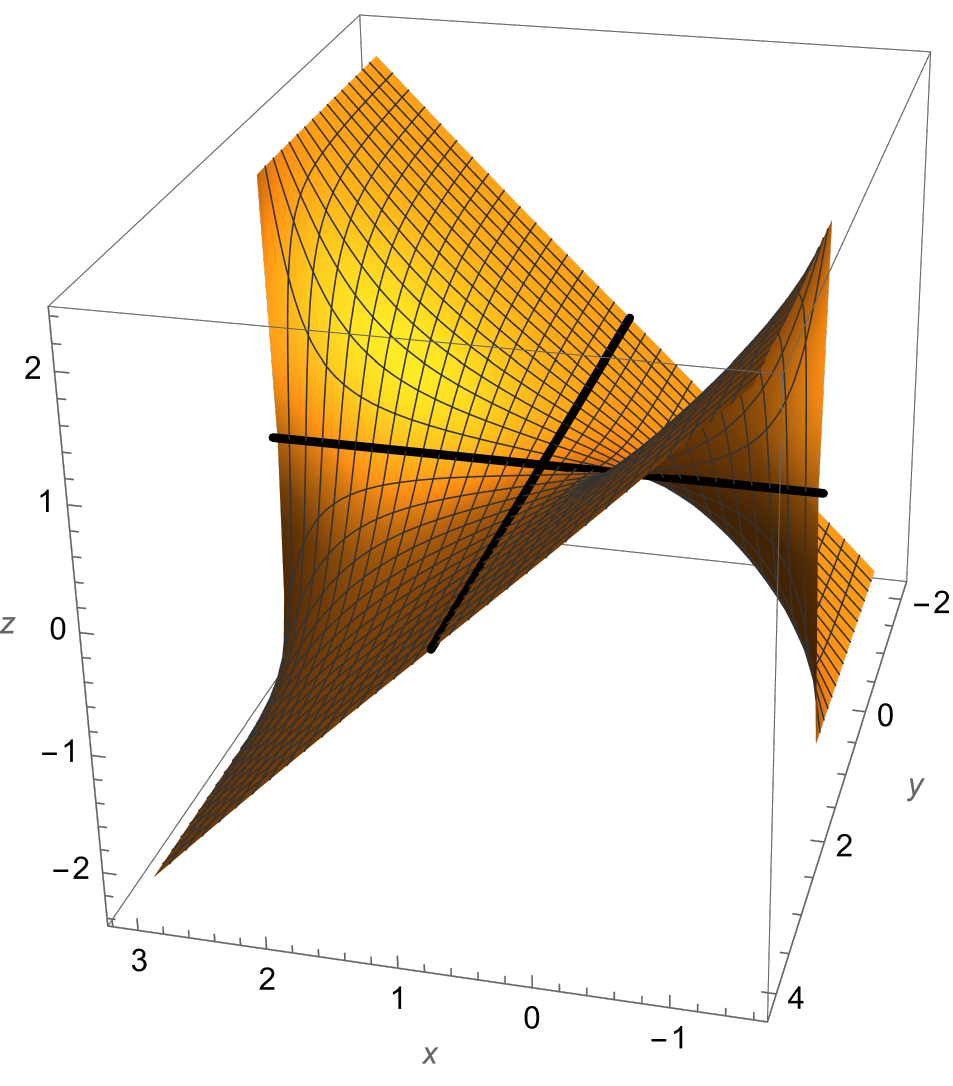}
 \end{center}
\caption{ The surface $\Sigma_3$ defined in \eqref{ex3}. Left: a piece of the surface showing the horizontal lines $L_i$, $1\leq i\leq 3$ (black) contained in  $\Sigma_3$.    Right: the same surface after reflections across horizontal lines.  }\label{fig6}
\end{figure}

We consider Eqs. \eqref{ss2} and take $k=2$, $r=a=d_1=d_2=c_1=1$. By a direct integration, we obtain
\begin{equation*}
\left\{
\begin{split}
f(x)&=2 \tanh ^{-1}\left(\text{sech} (\frac{1}{2} ) \tan  (\frac{x}{\sqrt{2}} )+\tanh  (\frac{1}{2} )\right),\\
g(y)&=\log \left(\frac{e \coth  (\frac{y}{\sqrt{2}} )-\sqrt{e}}{\sqrt{e}-\coth  (\frac{y}{\sqrt{2}} )}\right),\\
h(z)&=\log \left(\tan  (\frac{z}{\sqrt{2}} )\right)-1.
\end{split}
\right.
\end{equation*}
We determine the domains of the functions. Since the argument of  $\tanh^{-1}$ must be between $-1$ and $1$, we deduce that the domain of $f$ is $(x_1,x_2)$ with $x_1\approx -1.45$ and $x_2\approx 0.77$. For the function $g$, the argument of the logarithm  must be positive, so $y\in (y_1,\infty)$ with $y_1\approx 0.995$. Finally, $z\in (0,z_2)$, with $z_2=\frac{\pi}{\sqrt{2}}$. Moreover, from the expression for  $h(z)$ and substituting  into  \eqref{s1}, we obtain
\begin{equation}\label{ex3}
\Sigma_3: z=z(x,y)=\sqrt{2} \tan ^{-1}\left(e^{-f(x)-g(y)+1}\right).
\end{equation}
Since  $\lim_{x\to x_1}f(x)= -\infty$, $ \lim_{x\to x_2}f(x)=\infty$ and $ \lim_{y\to y_1}g(y)=\infty$, we conclude that the horizontal lines 
\begin{equation*}
\begin{split}
L_1&=\{(x_1,y,z_2):0<y<\infty\},\\
L_2&=\{(x_2,y,0):0<y<\infty\},\\
L_3&=\{(x,y_1,0):x_1<x<x_2\},
\end{split}
\end{equation*}
are limit points of $\Sigma_3$ (Fig. \ref{fig6}, left). Thus, the reflections of $180^0$ across $L_i$ ($1\leq i\leq 3$)  extend analytically $\Sigma_3$ (Fig. \ref{fig6}, right). 
  
  \begin{proposition} If $\Sigma_3$ is the surface  defined in  \eqref{ex3}, then $\Sigma_2$ is doubly periodic. The surface is a global graph on $\r^2$ except for a discrete set points.
\end{proposition}

\section*{Data availibility}
We do not analyze or generate any datasets, because our work proceeds within a theoretical and mathematical approach.

The author has no conflict of interest to declare that are relevant to the content of this article.
  %%%%%%%%%%%%%%%%%%%%%%%%
\section*{Acknowledgements}
The author has been partially supported by MINECO/MICINN/FEDER grant no. PID2023-150727NB-I00, and by the ``Mar\'{\i}a de Maeztu'' Excellence Unit IMAG, reference CEX2020-001105- M, funded by MCINN/AEI/10.13039/ 501100011033/ CEX2020-001105-M.

\end{document}